\documentclass[12pt]{amsart}
\usepackage{amscd,amssymb}
\usepackage[graph,frame,poly,arc]{xy}  
\usepackage[plainpages,backref,urlcolor=blue]{hyperref}

\theoremstyle{plain}
\newtheorem{thm}[subsection]{Theorem}
\newtheorem{lem}[subsection]{Lemma}
\newtheorem{prop}[subsection]{Proposition}
\newtheorem{cor}[subsection]{Corollary}

\theoremstyle{definition}
\newtheorem{rk}[subsection]{Remark}
\newtheorem{definition}[subsection]{Definition}
\newtheorem{ex}[subsection]{Example}

\newtheorem{question}[subsection]{Question}

\numberwithin{equation}{section}
\newcommand{\D}{{\mathcal D}}

\newcommand{\Q}{\mathbb{Q}}

\newcommand{\C}{\mathbb{C}}
\newcommand{\PP}{\mathbb{P}}

\newcommand{\N}{\mathbb{N}}

\DeclareMathOperator{\rank}{rank}

\begin{document}
%\date{June 4, 2009}

\title [Mixed multiplicities, Hilbert polynomials and homaloidal surfaces]
{Mixed multiplicities, Hilbert polynomials and homaloidal surfaces}

\author[Alexandru Dimca]{Alexandru Dimca$^1$}
\address{Univ. Nice Sophia Antipolis, CNRS,  LJAD, UMR 7351, 06100 Nice, France. }
\email{dimca@unice.fr}

\author[Gabriel Sticlaru]{Gabriel Sticlaru}
\address{Faculty of Mathematics and Informatics,
Ovidius University,
Bd. Mamaia 124, 900527 Constanta,
Romania}
\email{gabrielsticlaru@yahoo.com }
\thanks{$^1$ Partially supported by Institut Universitaire de France.}



\begin{abstract} We investigate the relationship among several numerical invariants associated to a (free) projective hypersurface $V$: the sequence of mixed multiplicities of its Jacobian ideal, the Hilbert polynomial of its Milnor algebra, and the sequence of exponents when $V$ is free.
As a byproduct, we obtain explicit equations for some of the homaloidal surfaces in the projective 3-dimensional space constructed by C. Ciliberto, F. Russo and  A. Simis.

\end{abstract}
 
\maketitle

%\tableofcontents

\section{Introduction} 

Let $S=\C[x_0,...,x_n]$ be the graded polynomial ring in $n+1$ indeterminates with complex coefficients, and let $S_k$ denote the vector space of degree $k$ homogeneous polynomials. For a degree $d$, consider a polynomial $f \in S_d$, the corresponding Jacobian ideal $J_f$ generated by the partial derivatives $f_j$ of $f$ with respect to $x_j$ for $j=0,...,n$ and the graded Milnor (or Jacobian) algebra $M(f)=\oplus_k M(f)_k=S/J_f$. 
The  Hilbert polynomial $P(M(f))$ of the graded $S$-module $M(f)$ encodes information on the 
 projective hypersurface $V=V(f):f=0$  in $\PP^n$ and the associated singular subscheme $\Sigma=\Sigma(f)$ defined by the Jacobian ideal. For example, when $V$ has only isolated singularities, then the Hilbert polynomial $P(M(f))$ is a constant, equal to $\tau(V)$, the total Tjurina number of $V$, see \cite{CD}.

On the other hand, one can consider the mixed multiplicities of  $\bf m$  and $J_f$,
\begin{equation}
\label{mm1}
 \mu^i(f)=e_i({\bf m}|J_f),
\end{equation}
where $\bf m$ denotes the maximal homogeneous ideal $(x_0,...,x_n)$ in $S$, see \cite{Trung} where the behavior of these invariants under generic hyperplane sections was studied and the recent paper \cite{JH} for applications to the topology of projective hypersurfaces. As noted already by B. Teissier \cite{Te} in the case of an isolated hypersurface singularity $(Y,0)$ and then proved by N. V. Trung \cite{Trung} for the non-isolated case, these mixed multiplicities have a topological flavor, i.e. they are related to the Milnor numbers of various generic linear sections of $(Y,0)$. In the projective global case of a hyperplane arrangement, these multiplicities coincide with the Betti numbers of the arrangement complement, see June Huh's paper \cite{JH}. These fundamental results on the mixed multiplicities and projective hypersurfaces are recalled in Theorem \ref{thmJH}. Next we list a number of direct corollaries of this result.

In the third section we discuss the relation between the Hilbert polynomial $P(M(f))$ of the Milnor algebra $M(f)$ and the sequence of mixed multiplicities $\mu^i(f)$ for $i=0,1,...,n$. In the case of a free hyperplane arrangement,  Theorem \ref{thm1} states that several data are equivalent: the exponents of the arrangement, the sequence of mixed multiplicities, the Hilbert polynomial $P(M(f))$ and the Hilbert function $H(M(f))$. In Example \ref{ex1} we describe two non-free plane arrangements in $\PP^3$ having the same sequence of mixed multiplicities, but distinct Hilbert polynomials. Example \ref{ex2} exhibits two free surfaces in $\PP^3$, one a plane arrangement, the other an irreducible surface, having the same exponents, but distinct multiplicity sequences.

In the case of plane arrangements we discuss in Theorem \ref{thm2} a relation, suggested by the case of free arrangements, between the Hilbert polynomial $P(M(f))$
and the Betti numbers of the complement. While this relation fails in part in general, it seems that it holds for the generic plane arrangements, though they are not free as soon as there are more than 4 planes in the arrangement, see Example \ref{exGA} and Question \ref{qGA}.

A projective hypersurface $V:f=0$ is said to be homaloidal if the degree of the gradient mapping
$$\phi_f=grad(f):\PP^n \dasharrow  \PP^n, \  \   \    \ x \mapsto (f_0(x):...:f_n(x))$$
is equal to one, in other words if $\phi_f$ is a birational isomorphism. This notion has attracted a lot of interest in the recent years, see \cite{CRS}, \cite{Do}, \cite{FM}, \cite{JH2}.
I. Dolgachev showed in \cite{Do} that there are exactly three homaloidal plane curves, up to projective equivalence:  (i) a smooth conic, (ii) a triangle and (iii) the union of a conic and one of its tangents. In higher dimensions $n\geq 3$, the homaloidal hypersurfaces  must have non-isolated singularities when the degree $d\geq 3$, a fact conjectured in \cite{DPap} and proved by June Huh in \cite{JH2} in full generality after a number of partial results along the years. 

Dolgachev's result shows that no irreducible free curve is homaloidal, since such a curve has degree $d\geq 5$, see \cite{ST}, \cite{DStFD}. Note however that a smooth conic is a nearly free curve
as defined in \cite{DStNF} and a triangle, respectively a conic plus tangent,  is a free curve arrangement. 

In the last section we investigate the relation between freeness and homaloidal property in the case of irreducible surfaces. Using three series of irreducible surfaces $D_d$, $D'_d$ and $D''_d$ introduced in \cite{DStFS} in our study of free and nearly free surfaces and an additional series $D'''_d$ described here, we construct for  $d= 7,8$ or
$d\geq 10$ four homaloidal surfaces in $\PP^3$ of degree $d$  not projectively equivalent to each other.
For $d=4,5,6$  (resp. $d=9$),  we give two  (resp. three) such homaloidal surfaces, see Corollary \ref{corhol}.
The proof that our four families are homaloidal is completely elementary, the only subtle point being in finding the right equations, and this was suggested by our study of free surfaces.
It turns out surprisingly that our homaloidal surfaces  belong to the class of surfaces $Y(1,d-1;d-2)$ constructed by 
C. Ciliberto, F. Russo and A. Simis in \cite{CRS} using dual hypersurfaces of some rational scroll
surfaces, see Remark \ref{rkY}. The new contribution of our note here is twofold: we give explicit equations for such surfaces and we point out the relation with the free and the nearly free surfaces.

The computations of various minimal resolutions and mixed multiplicities given in this paper were made using two computer algebra systems, namely CoCoA \cite{Co} and Singular \cite{Sing}.
The corresponding codes are available on request, some of them being available in \cite{St}.

\section{Hilbert polynomials and the mixed multiplicities $\mu^i(f)$} 

The {\it Hilbert function} $H(M(f)): \N \to \N$ of the graded $S$-module $M(f)$ is defined as usual by
\begin{equation}
\label{Hfunc}
 H(M(f))(k)= \dim M(f)_k.
\end{equation}

It is known that there is a unique polynomial $P(M(f))(t) \in \Q[t]$, called the {\it Hilbert polynomial} of $M(f)$, and an integer $k_0\in \N$ such that
\begin{equation}
\label{Hpoly}
 H(M(f))(k)= P(M(f))(k)
\end{equation}
for all $k \geq k_0$. 

For the definition and the basic properties of the mixed multiplicities $\mu^i(f)$ of  $\bf m$  and $J_f$,
where $\bf m$ denotes the maximal homogeneous ideal $(x_0,...,x_n)$ in $S$, we refer to \cite{JH} and the references in that paper. Here we just say that they are the coefficients of the (homogeneous) Hilbert polynomial $H(R({\bf m}|J_f) \in \Q[a,b]$  of the standard bigraded algebra
$$R({\bf m}|J_f)=\sum_{(a,b) \in \N^2}{\bf m}^aJ_f^b/{\bf m}^{a+1}J_f^b.$$
Hence, one of the leading ideas of this note, namely comparing the two Hilbert polynomials
$P(M(f)$ and $R({\bf m}|J_f)$, is a very natural one.
In the following result we list a number of fundamental properties of the mixed multiplicities $\mu^i(f)$. The first two of them are due to N. V. Trung, see Proposition 4.1 in \cite{Trung}, and the last two properties are proved by J. Huh in \cite{JH},  relying on results by S. Papadima and the first author in \cite{DPap}.

\begin{thm}
\label{thmJH} Let $V=V(f):f=0$ be a projective hypersurface in $\PP^n$ and denote by $D(f)$ the complement $\PP^n \setminus V$. Then the following hold.

\begin{enumerate}

\item If $V$ is smooth of degree $d$, then $\mu^i(f)=(d-1)^i$ for $i=0,...,n$.

\item If $H \subset \PP^n$ is a generic hyperplane and we denote by $V':f'=0$ the corresponding hyperplane section $V'=V \cap H$ in $H=\PP^{n-1}$, then $\mu^i(f')=\mu^i(f)$ for $i=0,1,...,n-1$.

\item For any hypersurface $V$ one has the inequalities $b_i(D(f)) \leq \mu^i(f)$ for $i=0,...,n$
and the Euler number $E(D(f))$ of the complement $D(f)$ is given by
$$ E(D(f))=\sum_{i=0}^n(-1)^i\mu^i(f).$$
Moreover if $V$ is a hyperplane arrangement, then $b_i(D(f)) =\mu^i(f)$ for $i=0,...,n$.

\item The degree of the gradient rational map 
$$\phi_f=grad(f):\PP^n \dasharrow \PP^n, \  \   \    \ x \mapsto (f_0(x):...:f_n(x))$$
 is equal to the last mixed multiplicity $\mu^n(f)$ if $\phi_f$ is dominant.

\end{enumerate}

\end{thm}
Using this powerful result, it is easy to obtain the following consequences. For the first one, see also Example 12 in \cite{JH}.

\begin{cor}
\label{corH1}
Let $V:f=0$ be a degree $d$ hypersurface having only isolated singularities. Then the multiplicity sequence $\mu^*(f)=(\mu^0(f),\mu^1(f),...,\mu^n(f))$ is given by 
$$(1,d-1,..., (d-1)^{n-1}, (d-1)^n-\mu(V)),$$
 where $\mu(V)$ is the total Milnor number of the hypersurface $V$, i.e. the sum of the Milnor numbers of the all singularities of $V$.
In particular,  the sequence $\mu^*(f)$ carries the same information as the Hilbert polynomial $P(M(f))$ if and only if all the singularities of $V$ are weighted homogeneous.

\end{cor}

\proof The claim about the sequence $\mu^*(f)$ follows using (1), (2) in Theorem \ref{thmJH} to determine the first $n$ entries and then the formula for the Euler number in (3)  in Theorem \ref{thmJH} to find the $(n+1)$-th entry.

The last claim follows from the fact mentionned in the Introduction that $P(M(f))=\tau(V)$, the total Tjurina number of $V$, and K. Saito's Theorem in \cite{KS0} saying that for an isolated hypersurface singularity $(Y,0)$ one has $\mu(Y,0)=\tau(Y,0)$ if and only if $(Y,0)$ is weighted homogeneous (in some coordinate system).

\endproof

\begin{cor}
\label{corH2}
The equality $b_1(D(f))=\mu^1(f)$ holds if and only if $V:f=0$ is a hyperplane arrangement.

\end{cor}

\proof

The claim follows from the fact that $b_1(D(f))=r-1$, where $r$ is the number of irreducible components of $V$, see \cite{D1}. Indeed, it is known that $\mu^1(f)=d-1$ for any reduced hypersurface.

\endproof

\begin{cor}
\label{corH3}
Let $V:f=0$ be a degree $d$ surface in $\PP^3$. Then the multiplicity sequence $\mu^*(f)=(\mu^0(f),\mu^1(f),\mu^2(f),\mu^3(f))$ is given by 
$$(1,d-1, (d-1)^2-\mu(C),  1+(d-1)(d-2)-\mu(C)-E(D(f))   ),$$
 where $C$ is a generic plane section of the surface $V$.

\end{cor}

\proof The claim about the sequence $\mu^*(f)$ follows using (1), (2) in Theorem \ref{thmJH} and Corollary \ref{corH1} to determine the first three entries and then the formula for the Euler number in (3)  in Theorem \ref{thmJH} to find the fourth entry.
\endproof

\section{The case of hyperplane arrangements} 

Corollary \ref{corH1} implies that for a line arrangement the sequence $\mu^*(f)$ carries the same information as the Hilbert polynomial $P(M(f))$. The higher dimension situations are described in the following results.

\begin{thm}
\label{thm1} Let $V=V(f):f=0$ be a  free hyperplane arrangement in $\PP^n$ and denote by $D(f)$ the complement $\PP^n \setminus V$. Then each of the following data determines the others.

\begin{enumerate}

\item the exponents $d_1 \leq d_2 \leq ... \leq d_n$ of the free arrangement.

\item the sequence $b_i(D(f))$ of Betti numbers of the complement $D(f)$, encoded in the Poincar\'e polynomial
$$P(D(f);t)=\sum_{i=0,n}b_i(D(f))t^{n-i}.$$

\item  the sequence $\mu^*(f)$ of mixed multiplicities. 

\item the Hilbert function $H(M(f))$ and the degree $d$ of $f$.

\item the Hilbert  polynomial $P(M(f))$ and the degree $d$ of $f$.

\end{enumerate}

\end{thm}

\proof By Theorem \ref{thmJH} the sequences (2) and (3) coincide for any hyperplane arrangement. If the arrangement is free, then the exponents $d_i$'s and the Poincar\'e polynomial
are related by the formula
$$P(D(f);t)=(t+d_1)(t+d_2) \cdots (t+d_n),$$
see \cite{OT}, \cite{Yo}. It follows that the sequence (1) determines the sequence (2) and conversely.
The exponents determine the minimal resolution for the Milnor algebra $M(f)$, namely for an essential arrangement $V$ one has the resolution
\begin{equation}
\label{res}
 0 \to \oplus_{j=1,n} S(-d-d_j+1) \to S(-d+1)^{n+1} \to S.
\end{equation}
The case of non essential arrangements can be reduced to the case of essential ones, and this is left for the interested reader as an exercise.
The resolution \eqref{res} clearly determines the Hilbert function for $M(f)$, so (1) determines (4), since in addition it is known that $d_1+d_2+...+d_n=d-1$, see for instance \cite{DStFS}.
By definition, (4) determines (5). So it remains to show that the Hilbert polynomial $P(M(f))$ and the degree $d$ determine the exponents.

The resolution \eqref{res} implies the following formula for the Hilbert polynomial
$$P(M(f))(x)={x+n \choose n}-(n+1){x+n+1-d \choose n}+ \sum_{j=1,n}{x+n+1-d-d_j \choose n},$$
where ${x+n \choose n}$ means here the polynomial in $x$ given by the usual formula
$$\frac{(x+n)(x+n-1)\cdots (x+1)}{n!},$$
and similarly for the other binomial coefficients. It follows that
$$P(M(f))(x-n-1+d)={x-1+d \choose n}-(n+1){x \choose n}+ \sum_{j=1,n}{x-d_j \choose n}.$$
Let $Q(x)=\sum_{j=1,n}{x-d_j \choose n}$ and note that to show that the polynomial $P(M(f))(x)$ determines the exponents is reduced to show that the polynomial $Q(x)$ determines the exponents. A simple computation shows that the coefficient of $x^{n-m}$ in $Q(x)$ for $0<m \leq n$ has the following form
$$(-1)^m\frac{{n \choose m}}{n!}T_m(d_1,...,d_n)+ \sum_{j=0,m-1}a_{mj}T_j(d_1,...,d_n),$$
where $a_{mj} \in \Q$ are known constants independent of the $d_j$'s and $$T_j(d_1,...,d_n)=\sum_{i=1,n}d_i^j.$$
 This formula shows that the polynomial $Q(x)$ determines all the symmetric functions $T_j(d_1,...,d_n)$, known as the Newton sums, for $j=1,...,n$, and hence $Q(x)$ determines the exponents $d_1 \leq d_2 \leq ...\leq d_n$.
\endproof

\begin{ex}
\label{ex1}
There are examples of non free plane arrangements $V:f=0$ and $V':f'=0$ in $\PP^3$ having the same sequence $\mu^*(f)= \mu^*(f')$ of mixed multiplicities, but distinct $P(M(f)) \ne P(M(f'))$
Hilbert polynomials.
Indeed, consider the following two distinct realizations of the Pappus configuration $9_3$. The first one is the line arrangement
$$W: g=xyz(x-y)(y-z)(x-y-z)(2x+y+z)(2x+y-z)(-2x+5y-z)=0,$$
and the second one is the line arrangement  given by
$$W': g'=xyz(x+y)(x+3z)(y+z)(x+2y+z)(x+2y+3z)(4x+6y+6z)=0.$$
 Note that $W$ and $W'$ consist both of 9 lines, and have the same number of double and triple points, i.e. 9 double points and 9 triple points. It follows that the complements have the same Betti numbers, e.g. using Theorem \ref{thmJH} (3) and Corollary \ref{corH1}.
Consider now the non essential arrangements $V:f=0$ and $V':f'=0$ in $\PP^3$ defined by the same equations, but now viewed in the polynomial ring with 4 indeterminates. 
The complement $D(f)$ (resp. $D(f')$) is just the product $D(g) \times \C$ (resp. $D(g') \times \C$), hence they have again the same Betti numbers, and thus one has $\mu^*(f)=\mu^*(f')$.
On the other hand,  the Hilbert-Poincar\'e polynomial of $M(f)$ is $45t-189$, 
and the Hilbert-Poincar\'e polynomial of $M(f')$ is $45t-190$, as explained in \cite{DStFS}.

One can get an example of essential arrangements using the plane arrangements $V:f=wg(x,y,z)=0$ and $V':f'=wg'(x,y,z)=0$, i.e. the usual cone construction in hyperplane arrangement theory, \cite{OT}. Now the complement $D(f)$ (resp. $D(f')$) is just the product $D(g) \times \C^*$ (resp. $D(g') \times \C^*$), hence they have again the same Betti numbers, and thus one has again $\mu^*(f)=\mu^*(f')$.
On the other hand,  the Hilbert-Poincar\'e polynomial of $M(f)$ is $54t-261$, 
and the Hilbert-Poincar\'e polynomial of $M(f')$ is $54t-262$, as explained in \cite{DStFS}.

\end{ex}

\begin{ex}
\label{ex2} If we consider free surfaces  and not only free hyperplane arrangements, then there are the following two possibilities for a free surface in $\PP^3$ of degree $7$ with exponents $d_1=1, d_2=2, d_3=3$: 

(i) the plane arrangement $V:f=(x^2-y^2)(x^2-z^2)(y^2-z^2)w=0$ and 

(ii) the free surface 
$$V':f'=x^6z+y^7+x^5yw+x^4y^3=0,$$ 
which corresponds to the degree $d=7$ in the sequence of surfaces $D_d$ discussed below in Proposition \ref{propDd}. 
Then the Milnor algebras $M(f)$ and $M(f')$ have the same Hilbert functions, but a direct computation shows that the multiplicity sequences are distinct, namely 
$$\mu^*(f)=(1,6,11,6)$$
and
$$\mu^*(f')=(1,6,7,1).$$
It follows that the exponents of $D:f=0$ do not determine the mixed multiplicity sequence
$\mu^*(f)$ for a free surface $D$.

\end{ex}

\begin{rk}
\label{rkGA} Consider the case of a free arrangement $V:f=0$ with exponents $(d_1,d_2,d_3)$ in $\PP^3$. Then the Poincar\'e polynomial
$$P(D(f);t)=\sum_{i=0,3}b_i(D(f))t^{3-i}$$
is given by the product $(t+d_1)(t+d_2)(t+d_3)$, see \cite{OT} and hence the mixed multiplicities
$\mu^*(f)$ are given by the sequence
$$(1=s_0,d-1=s_1,s_2,s_3)$$
where $s_i$ is the $i$-th elementary symmetric function in the exponents $d_i$'s.
Let $P(M(f))(k)=ak+b$ be the corresponding Hilbert polynomial. Then Theorem 4.6 in \cite{DStFS} implies that
\begin{equation}
\label{a}
 a=s_1^2-s_2=(d-1)^2-s_2= (d-1)^2-b_2(D(f)) ,
\end{equation}
and
\begin{equation}
\label{b}
 b=\frac{1}{2}[(d-1)^3-(3d-7)a-s_3]=\frac{1}{2}[(d-1)^3-(3d-7)a-b_3(D(f))].
\end{equation}

\end{rk}
One may wonder if the above relations among $a,b, b_2(D(f))$ and $b_3(D(f))$ are valid beyond the class of free arrangements. The answer is given by the following.

\begin{thm}
\label{thm2} Let $V=V(f):f=0$ be a  plane arrangement in $\PP^3$ and denote by $D(f)$ the complement $\PP^3 \setminus V$. Then the following hold.

\begin{enumerate}

\item The formula \eqref{a} involving the leading coefficient of the Hilbert polynomial $P(M(f))$ holds.

\item There are examples of  plane arrangements $V:f=0$  in $\PP^3$ for which the
 formula \eqref{b} involving the free term of the Hilbert polynomial $P(M(f))$ fails.

\end{enumerate}

\end{thm}

\proof

Note that the singular locus $\Sigma$ of $V$ is the union of all 1-dimensional edges $E$ which occur as intersections of planes in $V$. Such an edge $E$ is counted with a multiplicity $n_E=(k_E-1)^2$ where $k_E$ is the number of planes in $V$ passing through $E$. Since the leading coefficient $a$ is the degree of $\Sigma$, it follows that
$$a=\sum_En_E=\sum_E(k_E-1)^2.$$
On the other hand it follows from Theorem \ref{thmJH} that $b_2(D(f))=b_2(D(f'))$, where $V':f'=0$ is a generic plane section of $V$. Now $V'$ is a line arrangement and has as singularities one ordinary $k_E$-multiple point for each edge $E$ in $V$. Such a singularity has Milnor number
$\mu_E=(k_E-1)^2$ and Corollary \ref{corH1} implies that
$$b_2(D(f'))=(d-1)^2-\mu(V')=(d-1)^2-\sum_E\mu_E.$$
The previous formula for $a$ can be obtained also by recalling the general fact that $a$ is $P(M(f'))$ for any generic section.
This completes the proof of the first claim (1). To justify the claim (2) it is enough to consider again the plane arrangements discussed in Example \ref{ex1}. They have the same $b_2(D(f))$ but different $b$'s, so the formula cannot hold for both. A direct computation shows that in fact it holds for none.

\endproof

\begin{ex}
\label{exGA} Let $V:f=0$ be a generic plane arrangement in $\PP^3$, i.e. the intersection of any four  planes in $V$ is the empty set. We assume that $d$, the degree of $f$ is at least 4, and note that for $d=4$ we get the free arrangement $xyzw=0$, while for $d>4$ the corresponding arrangement is not free. This can be seen from the fact that the corresponding Poincar\'e polynomial
$$P(D(f);t)=t^3+(d-1)t^2+{d-1 \choose 2}t +{d \choose 3}-{d \choose 2}+{d \choose 1}-1,$$
computed for instance in \cite{OT}, is not a product of linear factors for $d>4$. Indeed, if it were the case, the obvious inequality $s_1^2 \geq 3s_2$ is equivalent to $(d-1)^2 \geq 3(d-1)(d-2)/3$, and this fails for $d>4$.
For this arrangement the formula \eqref{b} becomes the following.
\begin{equation}
\label{b1}
 b=\frac{1}{2}[(d-1)^3-(3d-7)(d-1)d/2-{d \choose 3}+{d \choose 2}-{d \choose 1}+1].
\end{equation}
This equality was checked by direct computation for generic plane arrangements of degree $d=5$,
e.g $f(x,y,z,w)=xyzw(x+y+z+w)$
and $d=6$, e.g. $ f(x,y,z,w)=xyzw(x+y+z+w)(x+2y+3z+4w)$.

\end{ex}
Motivated by this Example, we raise the following.

\begin{question}
\label{qGA}
Does the formula \eqref{b1} hold for any generic plane arrangement in $\PP^3$? In the affirmative case, is there a generalization to arbitrary dimension $n \geq 3$?
\end{question}
A positive answer would give a nice formula for the Hilbert polynomial $P(M(f))$ in this case.

\section{Free and nearly free surfaces and homaloidal surfaces}

In view of Theorem \ref{thmJH}, (4), a hypersurface $V:f=0$ in $\PP^n$ is homaloidal if and only if the gradient mapping $\phi_f$ is dominant and  $\mu^n(f)=1$. It seems to us that a mysterious relation exists between homaloidal surfaces and free and nearly free surfaces in $\PP^3$.
The simpliest instance of this relation  can be seen in following situation, perhaps well known to the specialists.

\begin{prop}
\label{propD3}
The discriminant $\D_3$ of the binary cubic form in $(u,v)$
$$P=xu^3+yu^2v+zuv^2+wv^3$$
is given by 
$$\D_3: \Delta_3=y^2z^2 - 4xz^3 - 4y^3w + 18xyzw - 27x^2w^2=0  $$
and is a free, homaloidal surface in $\PP^3$.

\end{prop}

\proof

The freeness of $\D_3$ is proved in \cite{DStFS}. To show that $\D_3$ is homaloidal, we can either check that $\mu^3(\Delta_3)=1$ using a computer algebra system, or use Corollary \ref{corH3}. Indeed,
$E(D(f))=E(\PP^3)-E(\D_3)=4-4=0$ since $\D_3$ is homeomorphic to $\PP^1\times \PP^1$, see
\cite{DStFS}. On the other hand, the generic plane section $C$ of $\D_3$ is a quartic curve with 3 cusps, as the singular set of $\D_3$ is just the twisted cubic curve in $\PP^3$ and the transversal singularity along this curve is a cusp. Hence $\mu(C)=6$ and this gives
$$\mu^3(f)=1+(d-1)(d-2)-\mu(C)+E(D(f))=1+6-6+0=1.$$

\endproof

Next we show that the free and nearly free surfaces in $\PP^3$ produce plenty of examples of homaloidal surfaces. The sequence of  surfaces $$D_d:f_d=x^{d-1}z+y^d+x^{d-2}yw+x^4y^{d-4}=0,$$ 
for $d \geq 4$ was introduced in Example 4.10 (i) in \cite{DStFS}. The surfaces in this sequence are free for $d=7,8$, and nearly free for all the other values $d \geq 4$, $d\ne 7,8$, see Proposition 5.6 in \cite{DStFS}.  A direct computation for $4 \leq d \leq 11$ shows that at least in this range one has the following.
\begin{equation}
\label{holo}
 \mu^*(f_d)=(1,d-1,d,1).
\end{equation}
In view of Theorem \ref{thmJH} (4) this suggests the following result.

\begin{prop}
\label{propDd}
The surface  $D_d:f_d=x^{d-1}z+y^d+x^{d-2}yw+x^4y^{d-4}=0$ is homaloidal for any $d \geq 4$.
\end{prop}

\proof

The proof is by direct computation. We have to show that for a generic point 
$p=(\alpha:\beta:\gamma:\delta)\in \PP^3$, the inverse image $\phi_f^{-1}(p)$ under the gradient mapping consists of exactly one point in $\PP^3$. We can assume $\gamma=1$ and then
the other coordinates $\alpha, \beta, \delta$ are uniquely determined by $p$. The equality $f_z=1$ implies that $x=\epsilon$, a $(d-1)$-th root of unity. The equation $f_w=\delta$ then implies $y=y_0\epsilon$, where $y_0=\delta$. Similarly, the equation $f_y=\beta$ yields $w=w_0\epsilon$, with $w_0=(\beta-d\delta^{d-1}-(d-4)\delta^{d-3}$.
Finally, the equation $f_x=\alpha$ yields $z=z_0\epsilon$, with
$$z_0=\frac{\alpha-(d-2)\delta w_0-4\delta^{d-4}}{d-1}.$$
It follows that the point $(x:y:z:w)$ in the inverse image $\phi_f^{-1}(p)$ is unique, since
$$(x:y:z:w)=(1:y_0:z_0:w_0).$$

\endproof

The sequence of  surfaces $D'_d:f_d=x^{d-1}z+y^d+x^{d-2}yw+x^{d-5}y^5=0$ for $d \geq 5$ was introduced in Example 4.10 (ii) and Proposition 4.11 in \cite{DStFS}. Such a surface is  nearly free  for $6 \leq d \leq 9$ and is free with exponents $(1,4,d-6)$ for $d \geq 10$. The computation of the mixed multiplicities $\mu^*(f_d)$ for
$6 \leq d \leq 13$ gives exactly the same values as in \eqref{holo}, which motivates the following result.

\begin{prop}
\label{propD'd}
The surface  $D'_d:f_d=x^{d-1}z+y^d+x^{d-2}yw+x^{d-5}y^5=0$ is homaloidal for any $d \geq 5$.
\end{prop}

\proof

The proof is exactly as for Proposition \ref{propDd} above.

\endproof

The sequence of nearly free surfaces $D''_d:f_d=x^{d-1}z+y^d+x^{d-2}yw=0$ for $d \geq 4$ was introduced in Proposition 5.5 in \cite{DStFS}. The corresponding exponents are $(1,1,d-2)$ and the computation of the mixed multiplicities $\mu^*(f_d)$ for
$4 \leq d \leq 13$ gives exactly the same values as in \eqref{holo}. We have the following result.

\begin{prop}
\label{propD"d}
The surface $D''_d:f_d=x^{d-1}z+y^d+x^{d-2}yw=0$ is homaloidal for any $d \geq 4$.
\end{prop}

\proof

The proof is exactly as for Proposition \ref{propDd} above. 

\endproof

The same method yields the following new series of surfaces, not considered in \cite{DStFS}.

\begin{prop}
\label{propD"'d}
The surface $D'''_d:f_d=x^{d-1}z+y^d+x^{d-2}yw+xy^{d-1}=0$ is nearly free with exponents $(1,2,d-3)$ and is homaloidal for any $d \geq 5$.
\end{prop}

\begin{rk}
\label{mu2} One can prove that $\mu^2(f)=d$ for all the defining equations of the surface $D$, where $D$ is one of the surfaces $D_d$, $D'_d$, $D''_d$ and $D'''_d$. One way to do this is to use the discussion in \cite{CRS}, p. 1794, where it is shown that the inverse of the gradient mapping $\phi_f$ is given by forms of degree $d$.
Another way consists to use the formula for the Euler number $E(D(f))$ given in Theorem \ref{thmJH}, (3), the fact that $\mu^3(f)=1$ and the fact that $E(D(f))=E(\PP^3) -E(D)=4-3=1.$
The equality $E(D)=3$ follows from a computation of Hodge-Deligne polynomials as explained in \cite{DStFS}, Proposition 6.3. 

On the other hand, $\mu^2(f)=d$ implies that $\mu(C)=(d-1)(d-2)-1$, i.e. a generic plane section of such surfaces $D$ is never a rational cuspidal curve.

\end{rk}

We remark that there are some overlaps in our three families: up to projective equivalence one has
$D_4=D''_4$, $D_5=D'_5=D''_5$, $D_6=D''_6$ and $D'_6=D'''_6$, and $D_9=D'_9$. One can show that these are the only overlaps, essentialy just by looking at the exponents and using the fact the these exponents and the type of freeness are invariants under projective equivalence.

\begin{cor}
\label{corhol}
The following irreducible surfaces of degree $d$ in $\PP^3$ are homaloidal.
\begin{enumerate}

\item For $d=4$: the free discriminant surface $\D_3$ and  the nearly free surface $D_4=D''_4$.

\item For $d=5$: the nearly free surfaces $D_5=D'_5=D''_5$ and $D'''_5$.

\item For $d=6$: the nearly free surfaces $D_6=D''_6$ and $D'_6=D'''_6$.

\item For $d=7,8$: the free surface $D_d$, and the nearly free surfaces $D'_d$, $D''_d$, $D'''_d$.

\item For $d=9$: the nearly free surfaces $D_9=D'_9$, $D''_9$ and $D'''_9$.

\item For $d\geq 10$:  the free surface $D'_d$ and  the nearly free surfaces $D_d$, $D''_d$ and $D'''_d$.

\end{enumerate}
\end{cor}

\begin{rk}
\label{rkY} If $D$ is one of the surfaces $D_d$, $D'_d$, $D''_d$ and $D'''_d$, then its singular locus is the line $L:x=y=0$. If $P_t$ denotes the plane through $L$ with equation $y-tx=0$, then the intersection $D \cap P_t$ consists of the line $L$, the directix, with multiplicity $(d-1)$ and a simple line denoted by $L_t$. It follows that $D$ is a surface of type $Y(1,d-1)$ in the notation of Section 3.2 in \cite{CRS}. Moreover, the blow-up $\overline D$ of $D$ along the line $L$ is smooth,
and the exceptional divisor consists of two lines, one with multiplicity $d-2$. Following the notation on p. 1794 in \cite{CRS}, it follows that $D=Y(1,d-1;d-2)$, the most degenerate cases of surfaces of type $Y(1,d-1) \subset \PP^3$.

\end{rk}

\begin{rk}
\label{rkhol} The nearly free surface $\D'_4:f'=0$ of degree 6 described in Proposition 5.7 in \cite{DStFS} as a special hyperplane section of the discriminant $\D_4$ in $\PP^4$ of the binary quartic forms in $(u,v)$ is not homaloidal, since a direct computation shows that $\mu^3(f')=6$. 
On the other hand, the cubic surface
$$ C: z^3-2yzw+xw^2=0$$
is homaloidal since $C=Y(1,2)$, see \cite{CRS}, p. 1791, but $C$ is neither free, nor nearly free.

\end{rk}

\end{document}